\documentclass[12pt]{article}
\usepackage[left=2.5cm,right=2cm,
    top=2cm,bottom=2cm,bindingoffset=0cm]{geometry}
\usepackage{latexsym,amssymb,amsmath}
\usepackage[dvips]{graphicx,color}
\usepackage{mathrsfs}

\usepackage{hyperref}

\usepackage{xcolor}
\def\blue{\color{blue}}
\def\red{\color{red}}

\date{}

\title{Comment on the preprint ArXiv:
https://doi.org/10.48550/arXiv.2501.17471 {\red Two-dimensional
Calderon Problem and Flat Metrics} by V.A.Sharafutdinov}
\author{M.I. Belishev}
\date{}

\begin{document}
\maketitle

\begin{abstract}
The notice provides some critical remarks on a paper
(preprint) by V.A.Sharafutdinov. The remarks mainly concern
Theorem 1.1 on determination of the Riemann surface from its
Dirichlet-to-Neuman map, which is declared as one of the basic
results of the paper.
\end{abstract}

\subsubsection*{About Theorem 1.1}

The citations are marked by color. References correspond to
numbering in the \cite{Shara}. Additional references are labeled
by [a1]-[a11].
\bigbreak

\noindent {\bf Page 1.}
\smallbreak

\noindent {\bf 1)} In Abstract, the author claims: {\red
``We prove the uniqueness theorem: A compact connected
two-dimensional Riemannian manifold $(M,g)$ with non-empty
boundary is determined by the data $(\partial M,
g_\partial,\Lambda_g)$ uniquely up to conformal equivalence.''}\, -- \,In
fact, this theorem is proved in \cite{B}, whereas the \cite{Shara}
contributes no new ideas and/or techniques to the proof.
\smallbreak

\noindent {\bf 2)}\,\,{\red ``Nevertheless, I am not completely
satisfied with Lassas--Uhlmann’s proof. Indeed, the analytic
continuation along curves is used in the proof. If the surface $M$
is not simply connected, the continuation can give a multi-valued
function. The difficulty is not discussed in \cite{LU}.''}\, -- \,
To be satisfied, read papers attentively. The proof of \cite{LU}
uses the real-analytic continuation of Green's functions, whose
singularities are located in the boundary collar outside the
manifold. Such functions are a priori single-valued; it is
basically impossible to obtain multi-valued functions in this way.
So, everything is correct in \cite{LU}. \smallbreak

\noindent {\bf Page 3.\,}
\smallbreak

\noindent {\bf 3)}\,\,{\red ``Belishev's proof is based on the
well known Gelfand theorem: a compact topological space $X$ is
uniquely, up to a homeomorphism, determined by the Banach algebra
$C(X)$ of continuous functions.''}\, -- \,Right, it is well known
but the author is not fully versed in this subject. The mentioned
general theorem provides not much for the surface determination,
whereas the key fact is that the holomorphic continuous function
algebra is {\it generic} (see \cite{B}, sec 0.4). \smallbreak

\noindent {\bf 4)}\,\, {\red ``In my opinion, there are serious
gaps in Belishev's proof. First of all Theorem 1.1 is formulated
in the general case but is proved in \cite{B} only for surfaces
with connected boundaries. Some other gaps in Belishev's proof are
mentioned below.}\, -- \,Indeed, in \cite{B} there is an
inaccuracy in the statement: the proof deals with the case of a
single-component boundary, what should be specified in the
formulation. However, to extend it to the multi-component case is
an easy exercise, and to call it a {\red serious gap} is not quite
serious. If the author insists, there are the details (here and in
the subsequent we keep the notation of \cite{B}).

{\blue Instead of Lemma 3, \cite{B}, one uses its
``non-integrated'' version which is valid without the
connectedness assumption: {\it a smooth function $\eta=f+ih$ on
$\Gamma$ is an element of $\mathcal{A}(\Gamma)$ if and only if }
\begin{equation}
\label{dfgkjjkdf}
\partial_\gamma h=\Lambda f, \ \partial_\gamma f=-\Lambda h \qquad  on \quad \Gamma.
\end{equation}
The necessity is obvious since (\ref{dfgkjjkdf}) is just the
restriction on $\Gamma$ of the Cauchy-Riemann conditions for the
harmonic function $u^f+iu^h$. Let us prove the sufficiency. Let
$U$ be any neighborhood in $M$ whose boundary contains the segment
$l$ of $\Gamma$. Due to to the Poincare lemma, there is the smooth
(up to $l$) function $V$ obeying $dV=\star du^f$ in $U$; then
(\ref{dfgkjjkdf}) implies $dV=du^h$ on $l$ and one can chose $V$
is such a way that $V=h$ on $l$. Then $V=u^h$ on $U$ due to the
uniqueness of solutions to the Cauchy problem for the Laplace
equation. In particular, we have $du^h=dV=\star du^f$ on $U$ and,
since $U$ is arbitrary, on the whole $M$. Thus, the sufficiency if
proved.}

Lemma 3, \cite{B} is obtained from the above statement by the
additional (and unnecessary) integration of equations
(\ref{dfgkjjkdf}) provided that $\Gamma$ is connected; the
connectedness assumption is not used anywhere else in \cite{B}. If
one prefers to deal with the Hilbert transform but avoid
difficulties with integration over multi-connected boundary, one
can just use $\partial_\gamma\Lambda^{-1}$ instead of
$\Lambda J$.
\smallskip

At last, to recover $\Omega$ with a multi-component
$\Gamma$, one can use the technique of the paper \cite{BBK_IP},
which provides the determination of the surface with holes with
isolated borders.
\smallbreak

\noindent {\bf 5)}\,\,{\red ``A similar result is obtained by
Henkin-Michel \cite{HM}. Roughly speaking, they consider metric
surfaces $(M, g)$ such that the boundary curve ... has singular
points like angles, but must be a real analytic curve between
singular points. The article is overloaded by technical details
related to singular points. It is not easy to understand the main
idea of their approach, if such an idea exists.''}\, -- \,The
question arises: did the author communicate with V. Michel
regarding the content of \cite{HM}? Such a contact could help him
to understand the article and avoid embarrassment like the one
mentioned in {\bf 2)}. If not, then the author should have been
more restrained and written something like "Unfortunately, I could
not understand this article." \smallbreak

\noindent {\bf 6)}\,\, {\red ``In the present paper we combine
Belishev's approach \cite{B} with the approach of \cite{Velosip}
to get a complete rigorous (I hope) proof of Theorem 1.1 in the
general case...''}\, -- \,The formulation of the central Theorem
1, \cite{B} is literally true in the general case. With the
exception of the aforementioned minor correction {\blue{\bf 4)}},
approach \cite{B} does not need to be completed and/or combined
with anything else. \smallbreak

\noindent {\bf 7)}\,\, {\red ``The case of non-orientable surfaces
is reduced to the oriented case in \cite{BK_JIIPP}.''}\, -- \,
This reference is incorrect: the case of general non-orientable
surfaces is studied in \cite{BK_SIAM} while \cite{BK_JIIPP} deals
with the determination of the Mobius strip only. There are
articles available on this topic that the author should have been
aware of. \smallbreak

\noindent {\bf 8)}\,\, {\red ``In Section 2 we prove that topology
of a compact metric surface is determined by the DN-map.''}\, --
\,This topology is determined by a single number (Euler
characteristic and/or Betti number). How to find the number from
the DN-map is shown in \cite{B}. In our opinion, everything
written on pages 4-9 (from {\red Let $(M,g)$ be a compact ...} up
to the end of sec 2) is unnecessary. It looks like an
attempt to impress the reader by "high Math". There is a suitable
Russian saying "to shoot a cannon at sparrows". Dealing with a
two-dimensional problem, just a reference on \cite{BS} would be
quite enough. \smallbreak

\noindent {\bf Page 9.\,}
\smallbreak

\noindent {\bf 9)}\,\,\,Used in the proof of {\red Proposition
3.1}, the equation {\red $K_\varphi=
e^{-2\varphi}(K-\Delta_g\varphi)$} is not an invention of the
author, but a two-dimensional version of the well-known Yamabe
equation (in 2d case the scalar curvature is proportional to the
Gaussian one). This point should be noted and provided with an
appropriate reference. The author observes that, for the
two-dimensional manifolds with the boundary, the Yamabe problem
(for zero curvature) takes the form of the Dirichlet problem
\begin{equation}
\label{Ppp}
\Delta_g\varphi=K \text{ in } \Omega, \qquad \varphi=0 \text{ on } \Gamma
\end{equation}
and, therefore, one can take the solution to the Calderon problem
provided by \cite{B} and endow it with the flat metric in the same
conformal class by the use of the solution to (\ref{Ppp}). By the
way, such a trick is not new: Yamabe equation was already used in
inverse problems for determination of the conformal-euclidean
metric: see [a2], sec 1.3.

\smallbreak

\noindent {\bf Page 9.\,}
\smallbreak

\noindent {\bf 10)}\,\, {\red ``Belishev in \cite{B} discusses a
complex structure on a compact metric surface with non-empty
boundary, without giving a precise definition. In my opinion it is
his systematic mistake.}\, -- \,Declaring this as a `{\red
systematic mistake}' the author did not explain how it could
affect the correctness of the reasoning of \cite{B} (actually, it
cannot). Although there is the difference between various precise
(and well-known: see, e.g., \cite{Chirka 1,Chirka 2,Jur})
definitions, the complex structure (regardless of the choice of
definition) is (trivially!) determined by the following complex
atlas $\{U_k,w_k\}_k$ on $\Omega$, where $U_k$-s constitute the
open cover of $\Omega$, $w_k: \ U_k\to\mathbb{C}$ is an embedding
and $w_i\circ w_j^{-1}: \ {\rm int}\big(w_j(U_i\cap U_j)\big)\to
{\rm int}\big(w_i(U_i\cap U_j)\big)$ is a biholomorphism admitting
extension to a diffeomorphism between $w_j(U_i\cap U_j)$ and
$w_i(U_i\cap U_j)$. Such an atlas is recovered in \cite{B}. So,
any of definitions of complex/almost-complex/conformal structure
is completely consistent with the arguments of \cite{B}.
\smallbreak

\noindent {\bf Page 18.\,}
\smallbreak

\noindent {\bf 11)}\,\, {\red ``Belishev does not prove Theorem
3.3 in \cite{B}. There is only one sentence that can be related to
Lemma 8.1. Namely, we read on page 180 of \cite{B}: “Step 4.
Using functions $w\in\mathcal{A}(\Omega)$ as local homeomorphisms
$\Omega\to\mathbb{C}$ recover on $\Omega$ the complex
structure.” The sentence is hard understandable and proves
nothing.''}\, -- \,A way to endow the spectrum of the algebra
$\mathcal A(\Omega)$ with local holomorphic coordinates is shown
in \cite{B}. Such coordinates determine an atlas on $\Omega$. The
atlas determines a complex structure (that is the equivalence
class of the compatible atlases), and no more comment on are
required. The key fact is that the Gelfand transforms $\hat{\eta}$
of elements $\eta=w|_\Gamma$ of $\mathcal A(\Gamma)$ are related
to the holomorphic functions $w$ on $\Omega$ via
$\hat{\eta}=w\circ\beta^{-1}$, where $\beta$ is a homemorphism
from $\Omega$ onto the spectrum $\Omega'$ of $\mathcal A(\Gamma)$.
Therefore, to determine an atlas (and then the complex structure)
on $\Omega'$, it suffices to simply repeat the same operations as
with the original $\Omega$. Everything on $\Omega$ is transferred
to $\hat {\mathcal A}(\Omega)$ by the homeomorphism $\beta$ and
vice versa. For instance, the compatibility of holomorphic charts
on $\hat{\mathcal A}(\Omega)$ is ensured by compatibility of the
ones on $\Omega$. In a sense, such a trick is relevant for any
definition of the structure.

The mentioned details are quite trivial in Complex Analysis and
thus the conciseness of \cite{B} cannot be considered as a ``{\red
gap in Belishev's proof}''. This will not become true even if one
repeats it systematically. \smallbreak

\noindent {\bf Pages 9-25.\,}
\smallbreak

\noindent {\bf 12)}\,\,As is proved in \cite{LU} and \cite{B}, the
conformal class of the metric is determined by its DN-map. At the
same time, each class contains a unique flat metric that can be
recovered  using the Yamabe equation. In our opinion, this
essentially exhausts the content of the rest of the preprint. Does
it need so long presentation?

\subsubsection*{Conclusion}
The paper {\red \cite{Shara}} contains several high-profile, but
not quite responsible claims. As for determination of the surface
from its DN map, in fact, it adds nothing substantial to the
results \cite{B}. The author should realistically evaluate the
place and relevance of his activity against the background of the
real possibilities and achievements of the approach \cite{B}. For
help, see \cite{BK_CAOT,BK_stab_imm,BK_stab,K_stab}.

In our opinion, one should try to avoid the claims like {\red ``I am
not completely satisfied with Lassas--Uhlmann’s proof.''}, {\red
``It is not easy to understand the main idea ..., if such an idea
exists.''}, {\red ``In my opinion it is his systematic mistake.''}. They
look strange and doubtful in a mathematical text and may lead to
misunderstandings.

Also, we would advise the author to activate personal contacts.
Consulting with colleagues is more productive than wasting time on
reinventing the wheel and effort on unnecessary polemics.

\end{document}